\documentclass{manuscript}

\begin{document}
\title{The Kazdan--Warner equation on a finite graph with a Neumann boundary condition}
\author{%
    Wenyuan Ma$^{1}$\thanks{Corresponding author. E-mail: \href{mailto:wenyuanma1203@ruc.edu.cn}{\nolinkurl{wenyuanma1203@ruc.edu.cn}}}, Anfeng Yuan$^{2}$\\
    \affiliation{$^{1}$Department of Mathematics, Renmin University of China, Beijing 100872, China}\\
    \affiliation{$^{2}$Institute of Fundamental and Interdisciplinary Sciences, Beijing Union University,}\\
    \affiliation{Beijing 100101, China}
}
\maketitle

%% ── Abstract ─────────────────────────────────────────
% Abstract
{\noindent\small\textbf{Abstract}\quad
    We study the Kazdan--Warner equation $\Delta u = c - h \mathrm{e}^u$ on a finite connected domain $\Omega$ of a locally finite graph, with Neumann boundary data $\partial u / \partial\nu = g$. We reduce the boundary condition to homogeneous form and use variational methods and upper and lower solutions to study solvability.
    Let $c_g \coloneqq |\Omega|^{-1}\int_{\partial\Omega} g\, \mathrm{d}\mu$, and let $\psi$ solve $\Delta\psi = c_g$ with boundary data $\partial\psi / \partial\nu = g$ and zero mean on $\Omega$.
    When $c = c_g$ and $h \not\equiv 0$, a solution exists if and only if $h$ changes sign and $\int_\Omega h \mathrm{e}^\psi\, \mathrm{d}\mu < 0$.
    When $c > c_g$, a solution exists if and only if $h$ is positive somewhere in $\Omega$.
    When $c < c_g$, solvability requires $\int_\Omega h \mathrm{e}^\psi\, \mathrm{d}\mu < 0$. Under this condition there is a threshold $-\infty \le c_-(h) < c_g$ such that a solution exists for $c_-(h) < c < c_g$ and no solution exists for $c < c_-(h)$. If $h \le 0$ and $h \not\equiv 0$, a solution exists for every $c < c_g$.
}

\vspace{1ex}
{\noindent\small\textbf{Keywords}\quad
    Kazdan--Warner equation; Neumann boundary condition; variational method; upper and lower solutions.}

\vspace{1ex}
{\noindent\small\textbf{2020 MSC}\quad
    35R02; 35J61; 35A15.}

%% ── Sections ─────────────────────────────────────────
\section{Introduction}
\label{sec:introduction}

The Kazdan--Warner equation $\Delta u = c - h \mathrm{e}^u$ arises in the problem of prescribing Gaussian curvature on closed surfaces~\cite{Kazdan1974}. Related existence results on compact Riemannian surfaces include those of Yang and Zhu~\cite{Yang2018} and the critical, sign-changing case studied by Sun and Zhu~\cite{Sun2024}. Xu considers prescribed Gaussian and geodesic curvature on compact Riemann surfaces with boundary~\cite{Xu2024}.

On finite graphs, Grigor'yan, Lin, and Yang establish solvability results for this equation in the cases $c > 0$, $c = 0$, and $c < 0$~\cite{Grigoryan2016b}.
For $c < 0$, subsequent results address existence at the critical threshold~\cite{Ge2017a} and multiplicity for a related parameterized equation~\cite{Liu2020}.
Other developments include the $p$-Laplacian analogue for $p > 1$~\cite{Ge2020a}, degree-theoretic approaches~\cite{Sun2022,Liu2024e}, and equations involving spectral parameters~\cite{Li2023}.

Variational methods also yield existence results for Yamabe-type and other nonlinear equations on graphs~\cite{Grigoryan2016,Grigoryan2017}, and for mean-field equations on locally finite graphs~\cite{Lin2022}.
Further work treats finite-energy solutions of Kazdan--Warner equations on $\mathbb{Z}^2$~\cite{Chen2025} and the Kazdan--Warner equation on finite hypergraphs~\cite{Zhang2026b}.

For graph boundary value problems, Bensoussan and Menaldi develop variational and nonvariational formulations of Dirichlet, Neumann, and mixed problems~\cite{Bensoussan2005}. Hinz and Schwarz establish solvability and representation results for Neumann problems on finite and infinite graphs~\cite{Hinz2022a}.
Spectral studies provide Neumann Cheeger estimates~\cite{Hua2018} and isocapacitary bounds for Dirichlet, Neumann, and Steklov eigenvalues~\cite{Hua2024}.
In metric random walk spaces, including weighted graphs, nonlinear evolution problems have been studied with nonhomogeneous Neumann conditions~\cite{Mazon2020} and nonlinear, possibly dynamical, boundary conditions~\cite{Solera2023}.

Shi and Yu study the spectral effects of boundary conditions by comparing Dirichlet, Neumann, and Laplacian eigenvalues on finite weighted graphs~\cite{Shi2025}.
In their notation, $\Delta$ is nonnegative and the Neumann eigenvalue problem takes the form
\begin{equation*}
    \begin{cases}
        \Delta u = \lambda u, & \text{in } \Omega, \\[3pt]
        \displaystyle \frac{\partial u}{\partial\nu} = 0, & \text{on } \partial\Omega.
    \end{cases}
\end{equation*}

For a Kazdan--Warner boundary value problem, Pinamonti and Stefani consider the equation
\begin{equation*}
    \begin{cases}
        -\Delta_p u + \alpha \mathrm{e}^{\beta u} = f, & \text{in } \Omega, \\[3pt]
        u = \varphi, & \text{on } \partial\Omega,
    \end{cases}
\end{equation*}
on bounded domains of locally finite weighted graphs with nonempty interior and boundary.
For $1 \le p < \infty$ and nonnegative coefficients $\alpha$ and $\beta$, they prove that this Dirichlet problem has at most one solution~\cite{Pinamonti2022}.

We study the Kazdan--Warner equation with prescribed Neumann data,
\begin{equation}
    \begin{cases}
        \Delta u = c - h \mathrm{e}^u, & \text{in } \Omega, \\[3pt]
        \displaystyle \frac{\partial u}{\partial\nu} = g, & \text{on } \partial\Omega,
    \end{cases}
    \label{eq:neumann}
\end{equation}
where $c \in \mathbb R$, $h \colon \Omega \to \mathbb R$, and $g \colon \partial\Omega \to \mathbb R$.

We next introduce the graph notation. Let $G = (V, E)$ be a connected, locally finite graph. We consider a finite connected domain $\Omega \subset V$ and its boundary
\begin{equation*}
    \partial\Omega \coloneqq \{y \in V \setminus \Omega \colon xy \in E \text{ for some } x \in \Omega\}, \qquad \overline\Omega \coloneqq \Omega \cup \partial\Omega.
\end{equation*}
Write $x \sim y$ when $xy \in E$. Let the edge weight $w_{xy} = w_{yx} \ge 0$ satisfy $w_{xy} > 0$ if and only if $x \sim y$, and let $\mu \colon V \to (0, \infty)$ be a vertex measure. For a function $u \colon V \to \mathbb R$, we use the Laplacian
\begin{equation*}
    \Delta u(x) \coloneqq \frac{1}{\mu(x)} \sum_{y \sim x} w_{xy} \bigl(u(y) - u(x)\bigr).
\end{equation*}
The associated gradient form is
\begin{equation*}
    \nabla u \cdot \nabla v(x) \coloneqq \frac{1}{2 \mu(x)} \sum_{y \sim x} w_{xy} \bigl(u(y) - u(x)\bigr) \bigl(v(y) - v(x)\bigr),
\end{equation*}
and $|\nabla u|(x) \coloneqq \sqrt{\nabla u \cdot \nabla u(x)}$. For $A \subset V$ and $f \colon A \to \mathbb R$, write
\begin{equation*}
    \int_A f\, \mathrm{d}\mu \coloneqq \sum_{x \in A} \mu(x) f(x)
\end{equation*}
whenever the sum is defined, and set $|\Omega| \coloneqq \int_\Omega 1\, \mathrm{d}\mu$.

For a function $u \colon \overline\Omega \to \mathbb R$, define its outward normal derivative by
\begin{equation*}
    \frac{\partial u}{\partial\nu}(x) \coloneqq \frac{1}{\mu(x)} \sum_{y \in \Omega} w_{xy} \bigl(u(x) - u(y)\bigr), \qquad x \in \partial\Omega.
\end{equation*}
For functions $u, \varphi \colon \overline\Omega \to \mathbb R$, define the Dirichlet form associated with $\Omega$ by
\begin{equation}
    \begin{split}
        \mathcal E_\Omega(u, \varphi) & \coloneqq \frac{1}{2} \sum_{\substack{x, y \in \Omega \\ x \sim y}} w_{xy} \bigl(u(y) - u(x)\bigr) \bigl(\varphi(y) - \varphi(x)\bigr) \\
        & \quad + \sum_{\substack{x \in \Omega,\, y \in \partial\Omega \\ x \sim y}} w_{xy} \bigl(u(y) - u(x)\bigr) \bigl(\varphi(y) - \varphi(x)\bigr).
    \end{split}
    \label{eq:dirichlet-form}
\end{equation}
With these definitions, summation by parts gives the discrete Green identity
\begin{equation}
    \int_\Omega (\Delta u) \varphi\, \mathrm{d}\mu = -\mathcal E_\Omega(u, \varphi) + \int_{\partial\Omega} \frac{\partial u}{\partial\nu} \varphi\, \mathrm{d}\mu.
    \label{eq:green}
\end{equation}

The weak formulation of~\eqref{eq:neumann} is
\begin{equation*}
    \mathcal E_\Omega(u, \varphi) + c \int_\Omega \varphi\, \mathrm{d}\mu = \int_{\partial\Omega} g \varphi\, \mathrm{d}\mu + \int_\Omega h \mathrm{e}^u \varphi\, \mathrm{d}\mu
\end{equation*}
for every test function $\varphi \colon \overline\Omega \to \mathbb R$. The effect of the boundary condition is seen by taking $\varphi \equiv 1$ in~\eqref{eq:green}, which gives
\begin{equation}
    \int_\Omega \Delta u\, \mathrm{d}\mu = \int_{\partial\Omega} \frac{\partial u}{\partial\nu}\, \mathrm{d}\mu.
    \label{eq:divergence}
\end{equation}
Define
\begin{equation}
    c_g \coloneqq \frac{1}{|\Omega|} \int_{\partial\Omega} g\, \mathrm{d}\mu.
    \label{eq:cg}
\end{equation}

To remove the inhomogeneous boundary data $g$, let $\psi$ solve
\begin{equation*}
    \begin{cases}
        \Delta \psi = c_g, & \text{in } \Omega, \\[3pt]
        \displaystyle \frac{\partial\psi}{\partial\nu} = g, & \text{on } \partial\Omega,
    \end{cases}
\end{equation*}
with the normalization $\int_\Omega \psi\, \mathrm{d}\mu = 0$. The compatibility condition follows from~\eqref{eq:cg}, and existence follows from~\cite[Section~3.4]{Bensoussan2005}.

Then $u$ solves~\eqref{eq:neumann} if and only if $v \coloneqq u - \psi$ solves
\begin{equation}
    \begin{cases}
        \Delta v = \widetilde c - \widetilde h \mathrm{e}^v, & \text{in } \Omega, \\[3pt]
        \displaystyle \frac{\partial v}{\partial\nu} = 0, & \text{on } \partial\Omega,
    \end{cases}
    \label{eq:reduced}
\end{equation}
where
\begin{equation*}
    \widetilde c \coloneqq c - c_g, \qquad \widetilde h \coloneqq h \mathrm{e}^\psi.
\end{equation*}
Since $\mathrm{e}^\psi > 0$, $\widetilde h$ and $h$ have the same pointwise sign.

With the graph and boundary data $g$ fixed, we obtain the following results.

\begin{theorem}
    \label{thm:main}
    Let $h \not\equiv 0$. Then~\eqref{eq:neumann} with $c = c_g$ has a solution if and only if $h$ changes sign and
    \begin{equation*}
        \int_\Omega h \mathrm{e}^\psi\, \mathrm{d}\mu < 0.
    \end{equation*}
\end{theorem}

\begin{theorem}
    \label{thm:positive}
    Let $c > c_g$. Then~\eqref{eq:neumann} has a solution if and only if $h$ is positive somewhere in $\Omega$.
\end{theorem}

\begin{theorem}
    \label{thm:negative}
    Let $c < c_g$.
    \begin{enumerate}
        \item[(i)] If~\eqref{eq:neumann} has a solution, then
              $\displaystyle \int_\Omega \widetilde h\, \mathrm{d}\mu < 0$.
        \item[(ii)] If $\displaystyle \int_\Omega \widetilde h\, \mathrm{d}\mu < 0$, then there exists a constant $-\infty \le c_-(h) < c_g$ depending on $h$ such that~\eqref{eq:neumann} admits a solution for every $c_-(h) < c < c_g$, and no solution for any $c < c_-(h)$.
    \end{enumerate}
\end{theorem}

\begin{theorem}
    \label{thm:negative-all}
    If $h \le 0$ on $\Omega$ and $h \not\equiv 0$, then the threshold in Theorem~\ref{thm:negative} satisfies $c_-(h) = -\infty$. Consequently,~\eqref{eq:neumann} admits a solution for every $c < c_g$.
\end{theorem}

The remainder of this paper is organized as follows. Section~\ref{sec:preliminaries} records the compactness, Poincar\'e, and Trudinger--Moser estimates used below. Sections~\ref{sec:critical},~\ref{sec:positive}, and~\ref{sec:negative} prove the results in the critical, supercritical, and subcritical regimes, respectively.

\section{Preliminaries}
\label{sec:preliminaries}

We use the finiteness of $\overline\Omega$, the connectedness of $\Omega$, and the fact that each boundary vertex is adjacent to an interior vertex. The Dirichlet form $\mathcal E_\Omega$ defined in~\eqref{eq:dirichlet-form} is symmetric and bilinear, and $\mathcal E_\Omega(u, u) \ge 0$. When a constant depends on $\overline\Omega$, this includes dependence on the edge weights and vertex measure there.

Define the Sobolev space on $\overline\Omega$ by
\begin{equation*}
    W^{1,2}(\overline\Omega) \coloneqq \left\{\, u \in \mathbb R^{\overline\Omega} \colon \mathcal E_\Omega(u, u) + \int_{\overline\Omega} u^2\, \mathrm{d}\mu < +\infty \,\right\}
\end{equation*}
and equip it with the inner product
\begin{equation*}
    \langle u, v\rangle_{W^{1,2}(\overline\Omega)} \coloneqq \mathcal E_\Omega(u, v) + \int_{\overline\Omega} uv\, \mathrm{d}\mu.
\end{equation*}
The corresponding norm is
\begin{equation}
    \|u\|_{W^{1,2}(\overline\Omega)} \coloneqq \sqrt{\mathcal E_\Omega(u, u) + \int_{\overline\Omega} u^2\, \mathrm{d}\mu}.
    \label{eq:sobolev-norm}
\end{equation}
Since $\overline\Omega$ is finite, $W^{1,2}(\overline\Omega)$ is exactly the set of all real-valued functions on $\overline\Omega$ and is a finite-dimensional Hilbert space.
For any function $u$ defined on $\Omega$, write
\begin{equation*}
    \overline u \coloneqq \frac{1}{|\Omega|} \int_\Omega u\, \mathrm{d}\mu.
\end{equation*}

Following the arguments for finite graphs in~\cite{Grigoryan2016b}, we record the corresponding compactness and embedding results for $\overline\Omega$.

\begin{lemma}
    \label{lem:compact}
    Every bounded sequence $\{u_j\}$ in $W^{1,2}(\overline\Omega)$ has a subsequence converging in $W^{1,2}(\overline\Omega)$.
\end{lemma}

\begin{proof}
    Let $\{u_j\}$ be bounded in $W^{1,2}(\overline\Omega)$. By~\eqref{eq:sobolev-norm}, the sequence is bounded in $L^2(\overline\Omega)$. Since $\mu(x) > 0$ for every $x \in \overline\Omega$, each scalar sequence $\{u_j(x)\}$ is bounded. The set $\overline\Omega$ is finite, so a diagonal argument yields a subsequence, still denoted by $\{u_j\}$, and a function $u \colon \overline\Omega \to \mathbb R$ such that $u_j(x) \to u(x)$ for every $x \in \overline\Omega$.

    Both $\int_{\overline\Omega} (u_j - u)^2\, \mathrm{d}\mu$ and $\mathcal E_\Omega(u_j - u, u_j - u)$ are finite sums of terms converging to zero. Therefore
    \begin{equation*}
        \|u_j - u\|_{W^{1,2}(\overline\Omega)} \longrightarrow 0.
    \end{equation*}
    Since every function on the finite set $\overline\Omega$ belongs to $W^{1,2}(\overline\Omega)$, this proves the assertion.
\end{proof}

\begin{lemma}
    \label{lem:poincare}
    There exists a constant $C > 0$ depending only on $\overline\Omega$ such that every $u \in W^{1,2}(\overline\Omega)$ with $\overline u = 0$ satisfies
    \begin{equation*}
        \int_{\overline\Omega} u^2\, \mathrm{d}\mu \le C \mathcal E_\Omega(u, u).
    \end{equation*}
\end{lemma}

\begin{proof}
    Suppose, for contradiction, that no such $C$ exists. After normalization, there is a sequence $\{u_j\}$ such that
    \begin{equation*}
        \overline{u_j} = 0, \qquad \int_{\overline\Omega} u_j^2\, \mathrm{d}\mu = 1, \qquad \mathcal E_\Omega(u_j, u_j) \longrightarrow 0.
    \end{equation*}
    By~\eqref{eq:sobolev-norm}, the last two conditions imply that $\{u_j\}$ is bounded in $W^{1,2}(\overline\Omega)$. Lemma~\ref{lem:compact} therefore yields, after passing to a subsequence, a function $u_0 \in W^{1,2}(\overline\Omega)$ such that
    \begin{equation*}
        u_j \longrightarrow u_0 \quad \text{in } W^{1,2}(\overline\Omega).
    \end{equation*}
    Consequently,
    \begin{equation*}
        \overline{u_0} = 0, \qquad \int_{\overline\Omega} u_0^2\, \mathrm{d}\mu = 1, \qquad \mathcal E_\Omega(u_0, u_0) = 0.
    \end{equation*}
    Each term of $\mathcal E_\Omega(u_0, u_0)$ in~\eqref{eq:dirichlet-form} is nonnegative, with $w_{xy} > 0$ on each edge. Hence $u_0(y) = u_0(x)$ for every edge with both endpoints in $\Omega$ and for every edge joining $\Omega$ to $\partial\Omega$. Since $\Omega$ is connected and every vertex in $\partial\Omega$ is adjacent to a vertex in $\Omega$, the function $u_0$ is constant on $\overline\Omega$. Its zero average over $\Omega$ then implies $u_0 \equiv 0$, contradicting $\int_{\overline\Omega} u_0^2\, \mathrm{d}\mu = 1$.
\end{proof}

\begin{lemma}
    \label{lem:trudinger}
    For every $\beta \in \mathbb R$, there exists a constant $C > 0$ depending only on $\beta$ and $\overline\Omega$ such that
    \begin{equation*}
        \int_\Omega \mathrm{e}^{\beta v^2}\, \mathrm{d}\mu \le C
    \end{equation*}
    for every $v \in W^{1,2}(\overline\Omega)$ with $\mathcal E_\Omega(v, v) \le 1$ and $\overline v = 0$.
\end{lemma}

\begin{proof}
    If $\beta \le 0$, then $\mathrm{e}^{\beta v^2} \le 1$ and the assertion follows with $C = |\Omega|$. Assume henceforth that $\beta > 0$. By Lemma~\ref{lem:poincare}, there exists a constant $C_P > 0$ depending only on $\overline\Omega$ such that
    \begin{equation*}
        \int_{\overline\Omega} v^2\, \mathrm{d}\mu \le C_P \mathcal E_\Omega(v, v) \le C_P.
    \end{equation*}
    Set
    \begin{equation*}
        \mu_{\min} \coloneqq \min_{x \in \overline\Omega} \mu(x) > 0.
    \end{equation*}
    For every $x \in \overline\Omega$, we have
    \begin{equation*}
        \mu_{\min} v(x)^2 \le \mu(x) v(x)^2 \le \int_{\overline\Omega} v^2\, \mathrm{d}\mu \le C_P,
    \end{equation*}
    and hence $\|v\|_{L^\infty(\overline\Omega)}^2 \le C_P / \mu_{\min}$. It follows that
    \begin{equation*}
        \int_\Omega \mathrm{e}^{\beta v^2}\, \mathrm{d}\mu \le |\Omega| \mathrm{e}^{\frac{\beta C_P}{\mu_{\min}}}.
    \end{equation*}
    This proves the desired estimate.
\end{proof}

\section{The case \texorpdfstring{$c = c_g$}{c = c\_g}}
\label{sec:critical}

\begin{proof}[Proof of Theorem~\ref{thm:main}]
    \textbf{Necessity.} Let $u$ solve~\eqref{eq:neumann}. Then $v \coloneqq u - \psi$ solves~\eqref{eq:reduced} with $\widetilde c = 0$. Since $\widetilde h \not\equiv 0$, we have $\Delta v \not\equiv 0$, so $v$ is nonconstant. Multiplying the interior equation by $\mathrm{e}^{-v}$, integrating over $\Omega$, and applying~\eqref{eq:green} gives
    \begin{equation*}
        -\int_\Omega \widetilde h\, \mathrm{d}\mu = \int_\Omega \mathrm{e}^{-v} \Delta v\, \mathrm{d}\mu = -\mathcal E_\Omega(v, \mathrm{e}^{-v}) > 0,
    \end{equation*}
    where the boundary term vanishes because $\partial v / \partial\nu = 0$. The strict inequality follows from the strict decrease of $s \mapsto \mathrm{e}^{-s}$ and the connectedness of $\overline\Omega$. Thus
    \begin{equation*}
        \int_\Omega h \mathrm{e}^\psi\, \mathrm{d}\mu = \int_\Omega \widetilde h\, \mathrm{d}\mu < 0.
    \end{equation*}
    Integrating~\eqref{eq:neumann} and using~\eqref{eq:divergence} and~\eqref{eq:cg}, we obtain
    \begin{equation*}
        \int_\Omega h \mathrm{e}^u\, \mathrm{d}\mu = c |\Omega| - \int_{\partial\Omega} g\, \mathrm{d}\mu = 0.
    \end{equation*}
    Since $\mathrm{e}^u > 0$ and $h \not\equiv 0$, the latter identity forces $h$ to change sign.

    \textbf{Sufficiency.} Assume $h$ changes sign and $\int_\Omega h \mathrm{e}^\psi\, \mathrm{d}\mu < 0$. Define
    \begin{equation*}
        \mathcal M_{1} \coloneqq \left\{\, v \in W^{1,2}(\overline\Omega) \colon \int_\Omega \widetilde h \mathrm{e}^v\, \mathrm{d}\mu = 0,\; \int_\Omega v\, \mathrm{d}\mu = 0 \,\right\}.
    \end{equation*}

    We claim $\mathcal M_{1} \neq \varnothing$. Since $h$ changes sign, there exists $x_{1} \in \Omega$ with $h(x_{1}) > 0$. Define $v_{1}(x_{1}) \coloneqq \ell$ and $v_{1}(x) \coloneqq 0$ for $x \in \overline\Omega \setminus \{x_{1}\}$. Since $\widetilde h(x_{1}) > 0$, for all sufficiently large $\ell$ we have
    \begin{equation*}
        \int_\Omega \widetilde h \mathrm{e}^{v_{1}}\, \mathrm{d}\mu = (\mathrm{e}^\ell - 1) \mu(x_{1}) \widetilde h(x_{1}) + \int_\Omega \widetilde h\, \mathrm{d}\mu > 0.
    \end{equation*}
    Define the continuous function $\phi \colon [0, 1] \to \mathbb R$ by
    \begin{equation*}
        \phi(t) \coloneqq \int_\Omega \widetilde h \mathrm{e}^{t v_{1}}\, \mathrm{d}\mu.
    \end{equation*}
    Then $\phi(1) > 0$, while $\phi(0) = \int_\Omega \widetilde h\, \mathrm{d}\mu < 0$. Hence there exists $t_{0} \in (0, 1)$ with $\phi(t_{0}) = 0$. Take
    \begin{equation*}
        v^\ast \coloneqq t_{0} (v_{1} - \overline{v}_{1}).
    \end{equation*}
    Then $\int_\Omega v^\ast\, \mathrm{d}\mu = 0$, and
    \begin{equation*}
        \int_\Omega \widetilde h \mathrm{e}^{v^\ast}\, \mathrm{d}\mu = \mathrm{e}^{-t_{0} \overline{v}_{1}} \phi(t_{0}) = 0.
    \end{equation*}
    Thus $v^\ast \in \mathcal M_{1}$, proving the claim.

    Now consider the minimization problem
    \begin{equation*}
        J(v) \coloneqq \frac{1}{2} \mathcal E_\Omega(v, v), \qquad a \coloneqq \inf_{v \in \mathcal M_{1}} J(v).
    \end{equation*}
    Take a minimizing sequence $\{v_{n}\} \subset \mathcal M_{1}$ with $J(v_{n}) \to a$. Since $\int_\Omega v_{n}\, \mathrm{d}\mu = 0$, Lemma~\ref{lem:poincare} gives
    \begin{equation*}
        \int_{\overline\Omega} v_{n}^2\, \mathrm{d}\mu \le C \mathcal E_\Omega(v_{n}, v_{n}) = 2 C J(v_{n}).
    \end{equation*}
    Thus $\{v_{n}\}$ is bounded in $W^{1,2}(\overline\Omega)$. By Lemma~\ref{lem:compact}, after passing to a subsequence, we may assume that $v_{n} \to v_{\infty}$ in $W^{1,2}(\overline\Omega)$. Since $\overline\Omega$ is finite, this convergence is pointwise as well. Hence
    \begin{equation*}
        \int_\Omega v_{\infty}\, \mathrm{d}\mu = 0, \qquad \int_\Omega \widetilde h \mathrm{e}^{v_{\infty}}\, \mathrm{d}\mu = \lim_{n \to \infty} \int_\Omega \widetilde h \mathrm{e}^{v_{n}}\, \mathrm{d}\mu = 0,
    \end{equation*}
    so $v_{\infty} \in \mathcal M_{1}$. Moreover,
    \begin{equation*}
        J(v_{\infty}) = \lim_{n \to \infty} J(v_{n}) = a.
    \end{equation*}
    Thus $v_{\infty}$ is a minimizer of $J$ over $\mathcal M_{1}$.

    We derive the Euler--Lagrange equation for $v_{\infty}$. The differentials of the two constraint functionals are linearly independent at $v_{\infty}$. Indeed, if
    \begin{equation*}
        A \int_\Omega \widetilde h \mathrm{e}^{v_{\infty}} \phi\, \mathrm{d}\mu + B \int_\Omega \phi\, \mathrm{d}\mu = 0
    \end{equation*}
    for every $\phi \in W^{1,2}(\overline\Omega)$, then $A \widetilde h \mathrm{e}^{v_{\infty}} + B = 0$ pointwise on $\Omega$. Integrating over $\Omega$ and using the constraint gives $B = 0$, and then $A = 0$ because $\widetilde h \not\equiv 0$. The Lagrange multiplier method therefore gives constants $\lambda, \gamma \in \mathbb R$ such that for every $\phi \in W^{1,2}(\overline\Omega)$,
    \begin{equation*}
        \mathcal E_\Omega(v_{\infty}, \phi) = \lambda \int_\Omega \widetilde h \mathrm{e}^{v_{\infty}} \phi\, \mathrm{d}\mu + \gamma \int_\Omega \phi\, \mathrm{d}\mu.
    \end{equation*}
    By~\eqref{eq:green}, taking $\phi$ to be the indicator of $\{x\}$, first for $x \in \Omega$ and then for $x \in \partial\Omega$, gives
    \begin{equation}
        \Delta v_{\infty} = -\lambda \widetilde h \mathrm{e}^{v_{\infty}} - \gamma \quad \text{in } \Omega,
        \label{eq:euler-lagrange}
    \end{equation}
    together with $\partial v_{\infty} / \partial\nu = 0$ on $\partial\Omega$.
    Integrating~\eqref{eq:euler-lagrange} over $\Omega$ and using~\eqref{eq:divergence}, the homogeneous Neumann condition, and the constraint $\int_\Omega \widetilde h \mathrm{e}^{v_{\infty}}\, \mathrm{d}\mu = 0$ gives $\gamma = 0$.

    The only constant function with zero interior mean is zero, and $\int_\Omega \widetilde h\, \mathrm{d}\mu < 0$ excludes zero from $\mathcal M_{1}$. Thus $v_{\infty}$ is nonconstant. Multiplying~\eqref{eq:euler-lagrange} by $\mathrm{e}^{-v_{\infty}}$ and using~\eqref{eq:green} gives
    \begin{equation*}
        -\lambda \int_\Omega \widetilde h\, \mathrm{d}\mu = \int_\Omega \mathrm{e}^{-v_{\infty}} \Delta v_{\infty}\, \mathrm{d}\mu = -\mathcal E_\Omega(v_{\infty}, \mathrm{e}^{-v_{\infty}}) > 0.
    \end{equation*}
    It follows that $\lambda > 0$. Setting $\widehat v \coloneqq v_{\infty} + \log \lambda$, we obtain a solution of~\eqref{eq:reduced} with $\widetilde c = 0$. Hence $u \coloneqq \widehat v + \psi$ solves~\eqref{eq:neumann} with $c = c_g$.
\end{proof}

\section{The case \texorpdfstring{$c > c_g$}{c > c\_g}}
\label{sec:positive}

\begin{proof}[Proof of Theorem~\ref{thm:positive}]
    \textbf{Necessity.} Let $u$ solve~\eqref{eq:neumann}. Then $v \coloneqq u - \psi$ solves~\eqref{eq:reduced} with $\widetilde c > 0$. Integrating~\eqref{eq:reduced} over $\Omega$ and using~\eqref{eq:divergence} and the homogeneous Neumann condition gives
    \begin{equation*}
        \int_\Omega \widetilde h \mathrm{e}^v\, \mathrm{d}\mu = \widetilde c |\Omega| > 0.
    \end{equation*}
    Hence $\widetilde h$ is positive somewhere. Since $\widetilde h = h \mathrm{e}^\psi$ and $\mathrm{e}^\psi > 0$, $h$ is positive somewhere in $\Omega$.

    \textbf{Sufficiency.} Assume $h(x_{0}) > 0$ for some $x_{0} \in \Omega$. Define
    \begin{equation*}
        \mathcal M_{2} \coloneqq \left\{\, v \in W^{1,2}(\overline\Omega) \colon \int_\Omega \widetilde h \mathrm{e}^v\, \mathrm{d}\mu = \widetilde c |\Omega| \,\right\}.
    \end{equation*}
    We claim $\mathcal M_{2} \neq \varnothing$. Set
    \begin{equation*}
        v_{\ell}(x) \coloneqq
        \begin{cases}
            \ell, & x = x_{0}, \\
            0, & x \in \overline\Omega \setminus \{x_{0}\},
        \end{cases}
        \qquad \widetilde v_{\ell}(x) \coloneqq -\ell \quad (x \in \overline\Omega).
    \end{equation*}
    Since $\widetilde h(x_{0}) = h(x_{0}) \mathrm{e}^{\psi(x_{0})} > 0$, we have
    \begin{equation*}
        \int_\Omega \widetilde h \mathrm{e}^{v_{\ell}}\, \mathrm{d}\mu \to +\infty, \qquad \int_\Omega \widetilde h \mathrm{e}^{\widetilde v_{\ell}}\, \mathrm{d}\mu = \mathrm{e}^{-\ell} \int_\Omega \widetilde h\, \mathrm{d}\mu \to 0,
    \end{equation*}
    as $\ell \to +\infty$. Hence for $\ell$ sufficiently large,
    \begin{equation*}
        \int_\Omega \widetilde h \mathrm{e}^{v_{\ell}}\, \mathrm{d}\mu > \widetilde c |\Omega|, \qquad \int_\Omega \widetilde h \mathrm{e}^{\widetilde v_{\ell}}\, \mathrm{d}\mu < \widetilde c |\Omega|.
    \end{equation*}
    Define the continuous function $\phi \colon [0, 1] \to \mathbb R$ by
    \begin{equation*}
        \phi(t) \coloneqq \int_\Omega \widetilde h \mathrm{e}^{t v_{\ell} + (1 - t) \widetilde v_{\ell}}\, \mathrm{d}\mu.
    \end{equation*}
    Then $\phi(0) < \widetilde c |\Omega| < \phi(1)$, so there exists $t_{0} \in (0, 1)$ with $\phi(t_{0}) = \widetilde c |\Omega|$. Thus $\mathcal M_{2} \neq \varnothing$.

    For $\widetilde c > 0$, we minimize the functional
    \begin{equation*}
        J(v) \coloneqq \frac{1}{2} \mathcal E_\Omega(v, v) + \widetilde c \int_\Omega v\, \mathrm{d}\mu
    \end{equation*}
    over $\mathcal M_{2}$. Let $b \coloneqq \inf\limits_{v \in \mathcal M_{2}} J(v)$. For any $v \in \mathcal M_{2}$, set $w \coloneqq v - \overline v$, so that $\int_\Omega w\, \mathrm{d}\mu = 0$ and $\mathcal E_\Omega(w, w) = \mathcal E_\Omega(v, v)$. Then
    \begin{equation*}
        \int_\Omega \widetilde h \mathrm{e}^w\, \mathrm{d}\mu = \widetilde c |\Omega| \mathrm{e}^{-\overline v} > 0,
    \end{equation*}
    hence
    \begin{equation}
        J(v) = \frac{1}{2} \mathcal E_\Omega(v, v) - \widetilde c |\Omega| \log \int_\Omega \widetilde h \mathrm{e}^w\, \mathrm{d}\mu + \widetilde c |\Omega| \log(\widetilde c |\Omega|).
        \label{eq:J-form}
    \end{equation}
    Suppose first that $\mathcal E_\Omega(w, w) > 0$. Set $\widetilde w \coloneqq w / \sqrt{\mathcal E_\Omega(w, w)}$. Then $\int_\Omega \widetilde w\, \mathrm{d}\mu = 0$ and $\mathcal E_\Omega(\widetilde w, \widetilde w) = 1$. For $\varepsilon > 0$, the inequality $ab \le \varepsilon a^2 + b^2 / (4 \varepsilon)$ gives $w \le \varepsilon \mathcal E_\Omega(w, w) + \widetilde w^2 / (4 \varepsilon)$. Applying Lemma~\ref{lem:trudinger} to $\widetilde w$ with $\beta = 1 / (4 \varepsilon)$, we obtain
    \begin{equation}
        \int_\Omega \mathrm{e}^w\, \mathrm{d}\mu \le C_\varepsilon \mathrm{e}^{\varepsilon \mathcal E_\Omega(w, w)}.
        \label{eq:exp-w}
    \end{equation}
    If $\mathcal E_\Omega(w, w) = 0$, then $w \equiv 0$ by Lemma~\ref{lem:poincare}, and~\eqref{eq:exp-w} holds upon taking $C_\varepsilon \ge |\Omega|$. Consequently,
    \begin{equation*}
        0 < \int_\Omega \widetilde h \mathrm{e}^w\, \mathrm{d}\mu \le \|\widetilde h^{+}\|_{L^\infty(\Omega)} \int_\Omega \mathrm{e}^w\, \mathrm{d}\mu \le C_\varepsilon \|\widetilde h^{+}\|_{L^\infty(\Omega)} \mathrm{e}^{\varepsilon \mathcal E_\Omega(w, w)},
    \end{equation*}
    where $\widetilde h^{+}(x) \coloneqq \max\{\widetilde h(x), 0\}$.
    Substitution into~\eqref{eq:J-form} yields
    \begin{equation*}
        J(v) \ge \frac{1}{2} \mathcal E_\Omega(v, v) - \widetilde c |\Omega| \varepsilon \mathcal E_\Omega(v, v) - C_1.
    \end{equation*}
    Choosing $\varepsilon = 1 / (4 \widetilde c |\Omega|)$, we obtain
    \begin{equation}
        J(v) \ge \frac{1}{4} \mathcal E_\Omega(v, v) - C_1.
        \label{eq:lower-bound}
    \end{equation}
    Hence $J$ is bounded below on $\mathcal M_{2}$ and $b > -\infty$. Take a minimizing sequence $\{v_{k}\} \subset \mathcal M_{2}$ with $J(v_{k}) \to b$. By~\eqref{eq:lower-bound}, $\mathcal E_\Omega(v_{k}, v_{k})$ is bounded. Write $w_{k} \coloneqq v_{k} - \overline v_{k}$. From the definition of $J$,
    \begin{equation*}
        \overline v_{k} = \frac{1}{\widetilde c |\Omega|} J(v_{k}) - \frac{1}{2 \widetilde c |\Omega|} \mathcal E_\Omega(v_{k}, v_{k}),
    \end{equation*}
    so $\{\overline v_{k}\}$ is bounded. Since $\int_\Omega w_{k}\, \mathrm{d}\mu = 0$, Lemma~\ref{lem:poincare} gives
    \begin{equation*}
        \int_{\overline\Omega} w_{k}^2\, \mathrm{d}\mu \le C \mathcal E_\Omega(w_{k}, w_{k}) = C \mathcal E_\Omega(v_{k}, v_{k}).
    \end{equation*}
    Thus $\{w_{k}\}$, and consequently $\{v_{k}\}$, is bounded in $W^{1,2}(\overline\Omega)$. By Lemma~\ref{lem:compact}, after passing to a subsequence, we may assume that $v_{k} \to v_{\infty}$ in $W^{1,2}(\overline\Omega)$. Since $\overline\Omega$ is finite, this convergence is pointwise as well. Therefore
    \begin{equation*}
        \int_\Omega \widetilde h \mathrm{e}^{v_{k}}\, \mathrm{d}\mu \longrightarrow \int_\Omega \widetilde h \mathrm{e}^{v_{\infty}}\, \mathrm{d}\mu, \qquad \mathcal E_\Omega(v_{k}, v_{k}) \longrightarrow \mathcal E_\Omega(v_{\infty}, v_{\infty}).
    \end{equation*}
    It follows that $v_{\infty} \in \mathcal M_{2}$ and $J(v_{\infty}) = b$.

    We derive the Euler--Lagrange equation for $v_{\infty}$. Define $F(v) \coloneqq \int_\Omega \widetilde h \mathrm{e}^v\, \mathrm{d}\mu - \widetilde c |\Omega|$. Since
    \begin{equation*}
        DF(v_{\infty})[1] = \int_\Omega \widetilde h \mathrm{e}^{v_{\infty}}\, \mathrm{d}\mu = \widetilde c |\Omega| > 0,
    \end{equation*}
    the constraint is regular at $v_{\infty}$. The Lagrange multiplier method therefore gives $\lambda \in \mathbb R$ such that for every $\phi \in W^{1,2}(\overline\Omega)$,
    \begin{equation*}
        \mathcal E_\Omega(v_{\infty}, \phi) + \widetilde c \int_\Omega \phi\, \mathrm{d}\mu = \lambda \int_\Omega \widetilde h \mathrm{e}^{v_{\infty}} \phi\, \mathrm{d}\mu.
    \end{equation*}
    By~\eqref{eq:green}, taking $\phi$ to be the indicator of $\{x\}$, first for $x \in \Omega$ and then for $x \in \partial\Omega$, gives
    \begin{equation*}
        \Delta v_{\infty} = \widetilde c - \lambda \widetilde h \mathrm{e}^{v_{\infty}} \quad \text{in } \Omega, \qquad \frac{\partial v_{\infty}}{\partial\nu} = 0 \quad \text{on } \partial\Omega.
    \end{equation*}
    By~\eqref{eq:divergence} and the boundary condition, integration over $\Omega$ yields
    \begin{equation*}
        \widetilde c |\Omega| = \lambda \int_\Omega \widetilde h \mathrm{e}^{v_{\infty}}\, \mathrm{d}\mu = \lambda \widetilde c |\Omega|,
    \end{equation*}
    where the last equality uses $v_{\infty} \in \mathcal M_{2}$. Thus $\lambda = 1$, so $v_{\infty}$ solves~\eqref{eq:reduced} and $u \coloneqq v_{\infty} + \psi$ solves~\eqref{eq:neumann}.
\end{proof}

\section{The case \texorpdfstring{$c < c_g$}{c < c\_g}}
\label{sec:negative}

We treat this case by the method of upper and lower solutions combined with monotone iteration. We first define lower and upper solutions.

\begin{definition}
    A function $u_{-} \in W^{1,2}(\overline\Omega)$ is called a \emph{lower solution} of~\eqref{eq:neumann} if
    \begin{equation*}
        \begin{cases}
            \Delta u_{-} - c + h \mathrm{e}^{u_{-}} \ge 0, & \text{in } \Omega, \\[4pt]
            \displaystyle \frac{\partial u_{-}}{\partial\nu} \le g, & \text{on } \partial\Omega.
        \end{cases}
    \end{equation*}
    A function $u_{+} \in W^{1,2}(\overline\Omega)$ is called an \emph{upper solution} of~\eqref{eq:neumann} if
    \begin{equation*}
        \begin{cases}
            \Delta u_{+} - c + h \mathrm{e}^{u_{+}} \le 0, & \text{in } \Omega, \\[4pt]
            \displaystyle \frac{\partial u_{+}}{\partial\nu} \ge g, & \text{on } \partial\Omega.
        \end{cases}
    \end{equation*}
\end{definition}

We next prove the maximum principle and the existence of a solution.

\begin{lemma}
    \label{lem:max-principle}
    Let $k \colon \Omega \to (0, \infty)$. If $w \colon \overline\Omega \to \mathbb R$ satisfies
    \begin{equation*}
        \begin{cases}
            (\Delta - k) w \ge 0, & \text{in } \Omega, \\[4pt]
            \displaystyle \frac{\partial w}{\partial\nu} \le 0, & \text{on } \partial\Omega,
        \end{cases}
    \end{equation*}
    then $w \le 0$ on $\overline\Omega$.
\end{lemma}

\begin{proof}
    Suppose, for contradiction, that $M \coloneqq \max\limits_{\overline\Omega} w > 0$, and choose $x_{0} \in \overline\Omega$ with $w(x_{0}) = M$. If $x_{0} \in \Omega$, then $\Delta w(x_{0}) \le 0$ and
    \begin{equation*}
        (\Delta - k) w(x_{0}) = \Delta w(x_{0}) - k(x_{0}) M < 0,
    \end{equation*}
    contradicting the hypothesis. If $x_{0} \in \partial\Omega$, then
    \begin{equation*}
        \frac{\partial w}{\partial\nu}(x_{0}) = \frac{1}{\mu(x_{0})} \sum_{y \in \Omega,\, y \sim x_{0}} w_{x_{0} y} \bigl( w(x_{0}) - w(y) \bigr) \ge 0,
    \end{equation*}
    because $w(x_{0}) - w(y) \ge 0$ for all $y$. Since the normal derivative is also nonpositive at $x_{0}$, the sum vanishes. Thus $w(y) = M$ for every $y \in \Omega$ adjacent to $x_{0}$. Since every boundary vertex is adjacent to at least one interior vertex, there exists $y_{0} \in \Omega$ with $w(y_{0}) = M > 0$. The interior argument then yields a contradiction. Hence $M \le 0$.
\end{proof}

\begin{lemma}
    \label{lem:monotone}
    If~\eqref{eq:neumann} has a lower solution $u_{-}$ and an upper solution $u_{+}$ satisfying $u_{-} \le u_{+}$ on $\overline\Omega$, then it admits a solution $u$ satisfying $u_{-} \le u \le u_{+}$ on $\overline\Omega$.
\end{lemma}

\begin{proof}
    Set $f(x, s) \coloneqq c - h(x) \mathrm{e}^s$, $k_{1}(x) \coloneqq \max\{\, 1, \, -h(x) \,\}$, and $k(x) \coloneqq k_{1}(x) \mathrm{e}^{u_{+}(x)}$. Define $L w \coloneqq (\Delta - k) w$ in $\Omega$ and introduce the linear Neumann operator
    \begin{equation*}
        \mathscr L \colon W^{1,2}(\overline\Omega) \longrightarrow \mathbb R^{\Omega} \times \mathbb R^{\partial\Omega}, \qquad \mathscr L w \coloneqq \left( Lw, \frac{\partial w}{\partial\nu} \right).
    \end{equation*}
    If $\mathscr L w = 0$, applying Lemma~\ref{lem:max-principle} to both $w$ and $-w$ gives $w \equiv 0$. Thus $\ker \mathscr L = \{0\}$. The domain and codomain of $\mathscr L$ have the same finite dimension $\#\overline\Omega$, so $\mathscr L$ is invertible.

    Set $u_{0} \coloneqq u_{+}$. For each $j \ge 0$, let $u_{j + 1}$ be the unique solution of
    \begin{equation}
        \begin{cases}
            L u_{j + 1} = f(x, u_{j}) - k u_{j}, & \text{in } \Omega, \\[4pt]
            \displaystyle \frac{\partial u_{j + 1}}{\partial\nu} = g, & \text{on } \partial\Omega.
        \end{cases}
        \label{eq:iteration}
    \end{equation}

    We first show $u_{1} \le u_{0}$. Since $u_{0} = u_{+}$ is an upper solution,
    \begin{equation*}
        \begin{cases}
            L(u_{1} - u_{0}) = f(x, u_{0}) - \Delta u_{0} \ge 0, & \text{in } \Omega, \\[4pt]
            \displaystyle \frac{\partial(u_{1} - u_{0})}{\partial\nu} = g - \frac{\partial u_{+}}{\partial\nu} \le 0, & \text{on } \partial\Omega.
        \end{cases}
    \end{equation*}
    Lemma~\ref{lem:max-principle} gives $u_{1} \le u_{0}$.

    For each $x \in \Omega$ and $s \le u_{+}(x)$, we have
    \begin{equation}
        k(x) + h(x) \mathrm{e}^{s} \ge k_{1}(x) \bigl(\mathrm{e}^{u_{+}(x)} - \mathrm{e}^{s}\bigr) \ge 0,
        \label{eq:iteration-coefficient}
    \end{equation}
    since $h(x) \ge -k_{1}(x)$ and $k(x) = k_{1}(x) \mathrm{e}^{u_{+}(x)}$.
    Assume $u_{j} \le u_{j - 1} \le \cdots \le u_{+}$ for some $j \ge 1$. By~\eqref{eq:iteration} and the mean value theorem,
    \begin{equation*}
        L(u_{j + 1} - u_{j}) = \bigl(k + h \mathrm{e}^{\xi_{j}}\bigr) (u_{j - 1} - u_{j}) \ge 0 \quad \text{in } \Omega,
    \end{equation*}
    where $\xi_{j}(x)$ lies between $u_{j}(x)$ and $u_{j - 1}(x)$ and the inequality follows from~\eqref{eq:iteration-coefficient}. The boundary condition in~\eqref{eq:iteration} gives $\partial(u_{j + 1} - u_{j}) / \partial\nu = 0$ on $\partial\Omega$. Lemma~\ref{lem:max-principle} therefore gives $u_{j + 1} \le u_{j}$, completing the induction.

    Next, we prove $u_{-} \le u_{j}$ for all $j$. This holds for $j = 0$. Assume $u_{-} \le u_{j}$. At each $x \in \Omega$, the lower solution inequality and~\eqref{eq:iteration} give
    \begin{equation*}
        L(u_{-} - u_{j + 1}) \ge f(x, u_{-}) - f(x, u_{j}) + k(u_{j} - u_{-}) = \bigl(k + h \mathrm{e}^{\eta_{j}}\bigr) (u_{j} - u_{-}) \ge 0,
    \end{equation*}
    where the mean value theorem gives $\eta_{j}(x)$ between $u_{-}(x)$ and $u_{j}(x)$ and the last inequality follows from~\eqref{eq:iteration-coefficient}. On $\partial\Omega$ we have $\partial(u_{-} - u_{j + 1}) / \partial\nu = \partial u_{-} / \partial\nu - g \le 0$.
    Lemma~\ref{lem:max-principle} yields $u_{-} \le u_{j + 1}$, completing the induction. Hence
    \begin{equation*}
        u_{-} \le u_{j + 1} \le u_{j} \le \cdots \le u_{+} \quad \text{on } \overline\Omega.
    \end{equation*}

    For each $x \in \overline\Omega$, the sequence $\{u_{j}(x)\}$ is decreasing and bounded below by $u_{-}(x)$, hence converges pointwise. Since $\overline\Omega$ is finite, the convergence is uniform. Let $u(x) \coloneqq \lim\limits_{j \to \infty} u_{j}(x)$. Passing to the limit in~\eqref{eq:iteration} gives~\eqref{eq:neumann}, and $u_{-} \le u \le u_{+}$ on $\overline\Omega$.
\end{proof}

We now apply these lemmas to prove Theorem~\ref{thm:negative}.

\begin{proof}[Proof of Theorem~\ref{thm:negative}]
    (i) Let $u$ solve~\eqref{eq:neumann}. Then $v \coloneqq u - \psi$ solves~\eqref{eq:reduced} with $\widetilde c < 0$. Multiplying the equation for $v$ by $\mathrm{e}^{-v}$, integrating over $\Omega$, and applying~\eqref{eq:green} with $\partial v / \partial\nu = 0$ gives
    \begin{equation*}
        -\int_\Omega \widetilde h\, \mathrm{d}\mu = -\mathcal E_\Omega(v, \mathrm{e}^{-v}) - \widetilde c \int_\Omega \mathrm{e}^{-v}\, \mathrm{d}\mu > 0,
    \end{equation*}
    since $\mathcal E_\Omega(v, \mathrm{e}^{-v}) \le 0$ and $\widetilde c < 0$.

    (ii) Assume that $\overline{\widetilde h} < 0$. We apply Lemma~\ref{lem:monotone} to~\eqref{eq:reduced} with $c$, $h$, and $g$ replaced by $\widetilde c$, $\widetilde h$, and $0$, respectively. For each $\widetilde c < 0$, the constant $v_{-} \coloneqq -A$ satisfies the homogeneous Neumann condition and is a lower solution for all sufficiently large $A > 0$, since
    \begin{equation}
        -\widetilde c + \widetilde h \mathrm{e}^{-A} \longrightarrow -\widetilde c > 0 \quad \text{uniformly on } \Omega \quad (A \to +\infty).
        \label{eq:negative-lower}
    \end{equation}
    For any given upper solution $v_{+}$, we may increase $A$ so that $v_{-} \le v_{+}$ on $\overline\Omega$. To construct an upper solution, take $w$ satisfying $\int_\Omega w\, \mathrm{d}\mu = 0$ and
    \begin{equation}
        \begin{cases}
            \Delta w = \overline{\widetilde h} - \widetilde h, & \text{in } \Omega, \\[4pt]
            \displaystyle \frac{\partial w}{\partial\nu} = 0, & \text{on } \partial\Omega.
        \end{cases}
        \label{eq:negative-auxiliary}
    \end{equation}
    Such a $w$ exists because $\int_\Omega (\overline{\widetilde h} - \widetilde h)\, \mathrm{d}\mu = 0$. Since $\overline{\widetilde h} < 0$, we have $\max\limits_{x \in \Omega} |\widetilde h(x)| > 0$. Moreover, $\mathrm{e}^{a w} \to 1$ uniformly on the finite set $\Omega$ as $a \to 0^+$. We may therefore choose $a > 0$ such that
    \begin{equation*}
        |\mathrm{e}^{a w(x)} - 1| \le \frac{-\overline{\widetilde h}}{2 \max\limits_{y \in \Omega} |\widetilde h(y)|} \quad \text{for every } x \in \Omega.
    \end{equation*}
    Set $b \coloneqq \log a$, $\widetilde c_{0} \coloneqq a \overline{\widetilde h} / 2 < 0$, and $v_{+} \coloneqq a w + b$. Then $\partial v_{+} / \partial\nu = 0$ and
    \begin{equation*}
        \Delta v_{+} - \widetilde c_{0} + \widetilde h \mathrm{e}^{v_{+}} = \frac{a \overline{\widetilde h}}{2} + a \widetilde h (\mathrm{e}^{a w} - 1) \le 0.
    \end{equation*}
    Thus~\eqref{eq:reduced} has a solution at $\widetilde c_{0}$. A solution at any negative parameter is an upper solution at every larger negative parameter. Together with~\eqref{eq:negative-lower} and Lemma~\ref{lem:monotone}, this gives a solution at each such parameter. Let $\widetilde c_{-}(h)$ be the infimum of the negative parameters for which~\eqref{eq:reduced} has a solution. Then $-\infty \le \widetilde c_{-}(h) \le \widetilde c_{0} < 0$. The problem has a solution for every $\widetilde c_{-}(h) < \widetilde c < 0$ and no solution for any $\widetilde c < \widetilde c_{-}(h)$. Setting $c_{-}(h) \coloneqq \widetilde c_{-}(h) + c_g$ completes the proof.
\end{proof}

\begin{proof}[Proof of Theorem~\ref{thm:negative-all}]
    Since $h \le 0$ and $h \not\equiv 0$, we have $\widetilde h \le 0$ and $\overline{\widetilde h} < 0$. Let $w$ solve~\eqref{eq:negative-auxiliary} as in the proof of Theorem~\ref{thm:negative}. For any $\widetilde c < 0$, choose $a > 0$ with $a \overline{\widetilde h} < \widetilde c$ and $b \in \mathbb R$ such that $\mathrm{e}^{a w + b} \ge a$ on $\Omega$. Set $v_{+} \coloneqq a w + b$. Then $\partial v_{+} / \partial\nu = 0$ and
    \begin{equation*}
        \Delta v_{+} - \widetilde c + \widetilde h \mathrm{e}^{v_{+}} = a \overline{\widetilde h} - \widetilde c + \widetilde h (\mathrm{e}^{a w + b} - a) \le 0.
    \end{equation*}
    Thus $v_{+}$ is an upper solution. By~\eqref{eq:negative-lower}, we may choose $A > 0$ so large that $v_{-} \coloneqq -A$ is a lower solution and $v_{-} \le v_{+}$ on $\overline\Omega$. Lemma~\ref{lem:monotone} gives a solution of~\eqref{eq:reduced} for every $\widetilde c < 0$. Hence $c_{-}(h) = -\infty$.
\end{proof}

%% ── Acknowledgements ─────────────────────────────────
% \input{sections/acknowledgements}

%% ── Bibliography ──────────────────────────────────────
\bibliographystyle{abbrv}
\bibliography{refs}
\end{document}